\documentstyle[12pt]{article}
\begin{document}
\title{  An extremal problem on potentially $K_{p_{1},p_{2},...,p_{t}}$-graphic sequences
\thanks{  Project Supported by NNSF of China(10271105), NSF of Fujian,
Science and Technology Project of Fujian, Fujian Provincial
Training Foundation for "Bai-Quan-Wan Talents Engineering" ,
Project of Fujian Education Department and Project of Zhangzhou
Teachers College.}}
\author{{Chunhui Lai}\\
{\small Department of Mathematics}\\{\small Zhangzhou Teachers
College, Zhangzhou} \\{\small Fujian 363000,
 P. R. of CHINA.}\\{\small e-mail: zjlaichu@public.zzptt.fj.cn}}
\date{}
\maketitle
\begin{center}
\begin{minipage}{120mm}
\vskip 0.1in
\begin{center}{\bf Abstract}\end{center}
 {A sequence $S$ is potentially $K_{p_{1},p_{2},...,p_{t}}$ graphical if it has
a realization containing a $K_{p_{1},p_{2},...,p_{t}}$ as a subgraph, where
$K_{p_{1},p_{2},...,p_{t}}$ is a complete t-partite graph  with partition sizes
$p_{1},p_{2},...,p_{t} (p_{1}\geq p_{2}\geq ...\geq p_{t} \geq 1)$. Let $\sigma(K_{p_{1},p_{2},...,p_{t}}, n)$ denote the smallest degree sum
such that every $n$-term graphical sequence $S$ with
$\sigma(S)\geq \sigma(K_{p_{1},p_{2},...,p_{t}}, n)$ is potentially $K_{p_{1},p_{2},...,p_{t}}$
graphical.  In this paper, we prove that $\sigma (K_{p_{1},p_{2},...,p_{t}}, n)\geq 
2[((2p_{1}+2p_{2}+...+2p_{t}-p_{1}-p_{2}-...-p_{i}-2)n
-(p_{1}+p_{2}+...+p_{t}-p_{i})(p_{i}+p_{i+1}+...+p_{t}-1)+2)/2]$ for $n
\geq p_{1}+p_{2}+...+p_{t}, i=2,3,...,t.$  }\par
\par
\end{minipage}
\end{center}
 \par
 \section{Introduction}
\par

  If $S=(d_1,d_2,...,d_n)$ is a sequence of
non-negative integers, then it is called  graphical if there is a
simple graph $G$ of order $n$, whose degree sequence ($d(v_1 ),$
$d(v_2 ),$ $...,$ $d(v_n )$) is precisely $S$. If $G$ is such a
graph then $G$ is said to realize $S$ or be a realization of $S$.
A graphical sequence $S$ is potentially $H$ graphical if there is
a realization of $S$ containing $H$ as a subgraph, while $S$ is
forcibly $H$ graphical if every realization of $S$ contains $H$ as
a subgraph. Let $\sigma(S)=d(v_1 )+d(v_2 )+... +d(v_n ),$ and
$[x]$ denote the largest integer less than or equal to $x$. If $G$
and $G_1$ are graphs, then $G\cup G_1$ is the disjoint union of
$G$ and $G_1$. If $G = G_1$, we abbreviate $G\cup G_1$ as $2G$.
Let $K_k$, and $C_k$ denote a complete graph on $k$ vertices, and
a cycle on $k$ vertices, respectively. Let $K_{p_{1},p_{2},...,p_{t}}$ denote a
complete t-partite graph  with partition sizes $p_{1},p_{2},...,p_{t}
(p_{1}\geq p_{2}\geq ...\geq p_{t} \geq 1).$\par

Given a graph $H$, what is the maximum number of edges of a graph
with $n$ vertices not containing $H$ as a subgraph? This number is
denoted $ex(n,H)$, and is known as the Tur\'{a}n number. This
problem was proposed for $H = C_4$ by Erd\"os [3] in 1938 and in
general by Tur\'{a}n [11]. In terms of graphic sequences, the
number $2ex(n,H)+2$ is the minimum even integer $l$ such that
every $n$-term graphical sequence $S$ with $\sigma (S)\geq l $ is
forcibly $H$ graphical. Here we consider the following variant:
determine the minimum even integer $l$ such that every $n$-term
graphical sequence $S$ with $\sigma(S)\ge l$ is potentially $H$
graphical. We denote this minimum $l$ by $\sigma(H, n)$. Erd\"os,\
Jacobson and Lehel [4] showed that $\sigma(K_k, n)\ge
(k-2)(2n-k+1)+2$ and conjectured that equality holds. They proved
that if $S$ does not contain zero terms, this conjecture is true
for $k=3,\ n\ge 6$. The conjecture is confirmed in [5],[7],[8],[9]
and [10].
 \par
 Gould,\ Jacobson and
Lehel [5] also proved that  $\sigma(pK_2, n)=(p-1)(2n-2)+2$ for
$p\ge 2$; $\sigma(C_4, n)=2[{{3n-1}\over 2}]$ for $n\ge 4$. Yin
and Li [12] gave sufficient conditions for a graphic sequence
being potentially $K_{r,s}$-graphic, and determined
$\sigma(K_{r,r},n)$ for $r=3,4.$ Lai [6] proved that  $\sigma
(K_4-e, n)=2[{{3n-1}\over 2}]$ for $n\ge 7$.\  In this paper, we
prove that $\sigma (K_{p_{1},p_{2},...,p_{t}}, n)\geq 
2[((2p_{1}+2p_{2}+...+2p_{t}-p_{1}-p_{2}-...-p_{i}-2)n
-(p_{1}+p_{2}+...+p_{t}-p_{i})(p_{i}+p_{i+1}+...+p_{t}-1)+2)/2]$ for $n
\geq p_{1}+p_{2}+...+p_{t}, i=2,3,...,t.$  \par

\section{ Main results.} \par
{\bf  Theorem 1.} $\sigma (K_{p_{1},p_{2},...,p_{t}}, n)\geq 
2[((2p_{1}+2p_{2}+...+2p_{t}-p_{1}-p_{2}-...-p_{i}-2)n
-(p_{1}+p_{2}+...+p_{t}-p_{i})(p_{i}+p_{i+1}+...+p_{t}-1)+2)/2]$ for $n
\geq p_{1}+p_{2}+...+p_{t}, i=2,3,...,t.$ 
\par
{\bf Proof.} We first consider  $p_{1}+p_{2}+...+p_{i}-2p_{i}$ is even. 
If $n-p_{i}-p_{i+1}-...-p_{t}+1$ is even, let $n-p_{i}-p_{i+1}-...-p_{t}+1=2m, $ By
Theorem 11.5.9 of [2], $K_{2m}$ is the union of one 1-factor $M$
and $m-1$ spanning cycles $C_{1}^{1}, C_{2}^{1}, ...,
C_{m-1}^{1}.$ Let
$$H=C_{1}^{1}\bigcup C_{2}^{1}\bigcup...\bigcup
C_{\frac{p_{1}+p_{2}+...+p_{i}-2p_{i}}{2}}^{1}+K_{p_{i}+p_{i+1}+...+p_{t}-1}$$ 
Then $H$ is a realization of
$((n-1)^{p_{i}+p_{i+1}+...+p_{t}-1}, (p_{1}+p_{2}+...+p_{t}-p_{i}-1)^{n-p_{i}-p_{i+1}-...-p_{t}+1}).$
 Since $((p_{1}+p_{2}+...+p_{t}-p_{t})^{p_{t}}, (p_{1}+p_{2}+...+p_{t}-p_{t-1})^{p_{t-1}},
 ..., (p_{1}+p_{2}+...+p_{t}-p_{1})^{p_{1}})$ is the degree
sequence of $K_{p_{1},p_{2},...,p_{t}},$
 $((n-1)^{p_{i}+p_{i+1}+...+p_{t}-1}, (p_{1}+p_{2}+...+p_{t}-p_{i}-1)^{n-p_{i}-p_{i+1}-...-p_{t}+1})$ 
 is not potentially
$K_{p_{1},p_{2},...,p_{t}}$ graphic. Thus
$\sigma (K_{p_{1},p_{2},...,p_{t}}, n)\geq (p_{i}+p_{i+1}+...+p_{t}-1)(n-1) 
+ (p_{1}+p_{2}+...+p_{t}-p_{i}-1)(n-p_{i}-p_{i+1}-...-p_{t}+1)+ 2
= 2[((2p_{1}+2p_{2}+...+2p_{t}-p_{1}-p_{2}-...-p_{i}-2)n
-(p_{1}+p_{2}+...+p_{t}-p_{i})(p_{i}+p_{i+1}+...+p_{t}-1)+2)/2].$ 
Next, If $n-p_{i}-p_{i+1}-...-p_{t}+1$ is odd, let $n-p_{i}-p_{i+1}-...-p_{t}+1=2m+1, $ By
Theorem 11.5.9 of [2], $K_{2m+1}$ is the union of  $m$ spanning cycles $C_{1}^{1}, C_{2}^{1}, ...,
C_{m}^{1}.$ Let
$$H=C_{1}^{1}\bigcup C_{2}^{1}\bigcup...\bigcup
C_{\frac{p_{1}+p_{2}+...+p_{i}-2p_{i}}{2}}^{1}+K_{p_{i}+p_{i+1}+...+p_{t}-1}$$ 
Then $H$ is a realization of
$((n-1)^{p_{i}+p_{i+1}+...+p_{t}-1}, (p_{1}+p_{2}+...+p_{t}-p_{i}-1)^{n-p_{i}-p_{i+1}-...-p_{t}+1}),$
 and  we are done as before. This completes
the discussion for $p_{1}+p_{2}+...+p_{i}-2p_{i}$ is even.
\par
Now we consider  $p_{1}+p_{2}+...+p_{i}-2p_{i}$ is odd. 
If $n-p_{i}-p_{i+1}-...-p_{t}+1$ is even, let $n-p_{i}-p_{i+1}-...-p_{t}+1=2m, $ By
Theorem 11.5.9 of [2], $K_{2m}$ is the union of one 1-factor $M$
and $m-1$ spanning cycles $C_{1}^{1}, C_{2}^{1}, ...,
C_{m-1}^{1}.$ Let
$$H=M\bigcup C_{1}^{1}\bigcup C_{2}^{1}\bigcup...\bigcup
C_{\frac{p_{1}+p_{2}+...+p_{i}-2p_{i}-1}{2}}^{1}+K_{p_{i}+p_{i+1}+...+p_{t}-1}$$ 
Then $H$ is a realization of
$((n-1)^{p_{i}+p_{i+1}+...+p_{t}-1}, (p_{1}+p_{2}+...+p_{t}-p_{i}-1)^{n-p_{i}-p_{i+1}-...-p_{t}+1}),$
and  we are done as before. 
Next, If $n-p_{i}-p_{i+1}-...-p_{t}+1$ is odd, let $n-p_{i}-p_{i+1}-...-p_{t}+1=2m+1, $ By
Theorem 11.5.9 of [2], $K_{2m+1}$ is the union of  $m$ spanning cycles $C_{1}^{1}, C_{2}^{1}, ...,
C_{m}^{1}.$ Let
$$C_{1}^{1}=x_{1}x_{2}...x_{2m+1}x_{1}$$
$$H=(C_{1}^{1}\bigcup C_{2}^{1}\bigcup...\bigcup
C_{\frac{p_{1}+p_{2}+...+p_{i}-2p_{i}+1}{2}}^{1}+K_{p_{i}+p_{i+1}+...+p_{t}-1})$$
$$
-\{x_{1}x_{2}, x_{3}x_{4}, ...,
x_{2m-1}x_{2m}, x_{2m+1}x_{1}\} $$
Then $H$ is a realization of
$((n-1)^{p_{i}+p_{i+1}+...+p_{t}-1}, (p_{1}+p_{2}+...+p_{t}-p_{i}-1)^{n-p_{i}-p_{i+1}-...-p_{t}},$
$ 
(p_{1}+p_{2}+...+p_{t}-p_{i}-2)^{1}).$
It is easy to see that 
 $((n-1)^{p_{i}+p_{i+1}+...+p_{t}-1}, (p_{1}+p_{2}+...+p_{t}-p_{i}-1)^{n-p_{i}-p_{i+1}-...-p_{t}}, 
(p_{1}+p_{2}+...+p_{t}-p_{i}-2)^{1})$
 is not potentially
$K_{p_{1},p_{2},...,p_{t}}$ graphic. Thus
$\sigma (K_{p_{1},p_{2},...,p_{t}}, n)\geq (p_{i}+p_{i+1}+...+p_{t}-1)(n-1) 
+ (p_{1}+p_{2}+...+p_{t}-p_{i}-1)(n-p_{i}-p_{i+1}-...-p_{t})+(p_{1}+p_{2}+...+p_{t}-p_{i}-2)+ 2
= 2[((2p_{1}+2p_{2}+...+2p_{t}-p_{1}-p_{2}-...-p_{i}-2)n
-(p_{1}+p_{2}+...+p_{t}-p_{i})(p_{i}+p_{i+1}+...+p_{t}-1)+2)/2].$
 This completes
the discussion for $p_{1}+p_{2}+...+p_{i}-2p_{i}$  is  odd,
 and so finishes the proof of Theorem
1.
\par

   \par
       By  [4],[5],[7],[8],[9]
and [10]. The equality holds for $p_{1}=p_{2}=...=p_{t}=1,i=2$.

       \par

\end{document}